\documentclass[12pt]{article} 
\usepackage{amssymb, amsthm, amsmath, latexsym, verbatim, enumerate}
\usepackage[vcentermath, enableskew]{youngtab}
\theoremstyle{plain} 
\newtheorem*{conj}{Conjecture}
\newtheorem*{thm}{Theorem}

\newcommand{\cancel}{\hspace{-0.7em}\not}
\newcommand{\cneg}{\not\hspace{-0.2em}-}
\DeclareMathOperator{\lspec}{lspec}

\title{Sums of Consecutive Integers}
\author{Wai Yan Pong\\ California State University, Dominguez Hills}

\begin{document}
\maketitle

While looking for exercises for a number theory class, I recently came
across the following question in a book by Andr\'e Weil~\cite[Question
III.4]{weil}:
\begin{quote}
  {\em Which natural numbers can be written as the sum of two or more
    consecutive integers?}
\end{quote} 
The origin of this question is unknown to me, but one can easily
believe that it is part of the mathematical folklore. Solutions to
this problem and some generalizations maybe found in~\cite{guy}
and~\cite{lev}, for example. We will give a somewhat different proof
here, one that we hope readers will find intuitively appealing.

To get a feeling for the problem, let us consider some small values.
The list below gives one attempt at expressing the numbers up to 16 as
sums of two or more consecutive integers. We leave the right-hand side
of the equation blank if there is no such expression for the number on
the left-hand side.
\begin{align*}
  1& =\ ,&     2& =\ ,&         3& =1+2,&    4& =\ ,    \\
  5& =2+3,&    6& =1+2+3,&      7& = 3+4,&   8& =\ ,   \\
  9& =4+5,&   10& =1+2+3+4,&    11& =5+6,&   12& =3+4+5,\\
  13& =6+7,& 14& =2+3+4+5,& 15& =7+8,& 16& =\ ,
\end{align*}
A conjecture naturally arises from this list:
\begin{conj}
  A natural number is a sum of consecutive integers if and only
    if it is not a power of $2$.
\end{conj}
So, is it true? If so, what makes the powers of 2 special?

\section*{The Proof}
We start by defining a {\em decomposition of a natural number $ n$} to
be a sequence of consecutive natural numbers whose sum is $n$.  The
number of terms is called the {\em length of the decomposition}, and a
decomposition of length $1$ is called {\em trivial}.  Further, a
decomposition is called {\em odd} ({\em even}) if its length is odd
(even). For example, consider the number $15$. It has four
decompositions, $(15), (7, 8), (4, 5, 6)$ and $(1, 2, 3, 4, 5)$, one
even and three odd.

We will construct a one-to-one correspondence between the odd factors
of a number and its decompositions. This proves the conjecture since
the powers of $2$ are precisely those numbers with only one odd factor
(namely $1$) and thus they have only the trivial decomposition.

Here is our construction: Let $k$ be an odd factor of $n$. Then since
the sum of the $k$ integers from $-(k-1)/2$ to $(k-1)/2$ is $0$,
adding $n/k$ to each of them gives the sequence
\begin{equation} \label{eq:odd} 
\frac{n}{k}-\frac{k-1}{2},\quad \frac{n}{k}-\frac{k-1}{2}+1, \quad
\cdots, \quad \frac{n}{k}+\frac{k-1}{2},
\end{equation}
whose sum is 
\[ \sum_{j=-(k-1)/2}^{(k-1)/2} (\frac{n}{k}+j) 
= \frac{n}{k}\cdot k \quad + \sum_{j=-(k-1)/2}^{(k-1)/2} j = n+0 = n.\]
There are now two cases to consider:
\begin{itemize}
\item[(i)] $(k-1)/2 < n/k$. Then~\eqref{eq:odd} is already an odd
  decomposition of $n$. The length of this decomposition is $k$.

\item[(ii)] $(k-1)/2 \ge n/k$. Then~\eqref{eq:odd} begins with 0 or a
  negative number. After dropping the 0 and canceling the negative
  terms with the corresponding positive ones, we are left with the
  sequence
\begin{equation} \label{eq:even} 
  \frac{k-1}{2}- \frac{n}{k} +1, \quad \frac{k-1}{2}-\frac{n}{k} + 2,
  \quad \cdots, \quad \frac{n}{k} +
  \frac{k-1}{2},
\end{equation}
which is an even decomposition of $n$ of length $2n/k$.
\end{itemize} 
We call the decomposition in either~\eqref{eq:odd} or~\eqref{eq:even}
the {\em decomposition of $ n$ associated with $ k$}.

To show that every decomposition of $n$ arises this way, suppose
$(a+i)_{i=1}^m$ is a decomposition of $n$. Since its sum $m(2a+m+1)/2$
is $n$, either $m$ or $2a + m+1$ (but not both), depending on the
parity of $m$, is an odd factor of $n$.  It is then straightforward to
verify that the given decomposition is the one associated with that
odd factor.

We have in fact delivered more than we promised. Note that the
condition $(k-1)/2 < n/k$ is equivalent to $k-1 < 2n/k$, but since $k$
is odd the condition is therefore the same as $k < 2n/k$, or
equivalently $k < \sqrt{2n}$.  Likewise, $(k-1)/2 \ge n/k$ is
equivalent to $k > \sqrt{2n}$. Therefore, we have proved the
following:
\begin{thm} 
  There is a one-to-one correspondence between the odd factors of a
  natural number $n$ and its decompositions. More precisely, for each
  odd factor $k$ of $n$, if $k < \sqrt{2n}$, then~\eqref{eq:odd} is an
  odd decomposition of $n$ of length $k$. If $k > \sqrt{2n}$,
  then~\eqref{eq:even} is an even decomposition of $n$ of length
  $2n/k$. Moreover, these are all the decompositions of $n$.
\end{thm}
An immediate consequence of the theorem is that {\em the number of
  decompositions of $n$ is the number of odd factors of $n$}, which is
$\prod_p (e_p+1)$, where $p$ runs through the odd primes and $e_p$ is
the power of $p$ appearing in the prime factorization of $n$. Let us
illustrate both this formula and the theorem with the example $n=45$.
Since $45=3^2 \cdot 5$, there should be six decompositions, four odd
and two even, and indeed here they are:
\begin{table}[h]
\caption{\small The decompositions of 45} \label{t:decomp45}
\begin{center}
\begin{tabular}{c c c c c c}
\hline
$k$ & $(k-1)/2$ & $n/k$ & decomposition & length & parity\\
\hline
1 & 0 & 45 & (45) & 1 &odd\\
3 & 1 & 15 & (14, 15, 16) & 3 &odd\\
5 & 2 & 9  & (7, 8, 9, 10, 11) & 5 &odd\\
9 & 4 & 5  & (1, 2, 3, 4, 5, 6, 7, 8, 9) & 9 &odd\\
15& 7 & 3  & (5, 6, 7, 8, 9, 10) &6 &even\\
45& 22& 1  & (22, 23) & 2 &even\\
\hline
\end{tabular}
\end{center}
\end{table}

\section*{The Length Spectra}
We define the {\em length spectrum of $ n$}, denoted by
$\lspec(n)$, to be the set of lengths of the decompositions of $n$.
According to the theorem, $\lspec(n)$ is the set
\begin{equation*} 
\left\{ k \colon k \ \text{odd},\ k\, |\, n,\ k < \sqrt{2n} \right\} 
\cup \left\{ 2n/k \colon k\ \text{odd},\ k\, |\, n,\ 
k > \sqrt{2n} \right\}.
\end{equation*}
For example, we have seen that $\lspec(45) = \{1,3,5,9\} \cup \{2,6\}$
(Table~\ref{t:decomp45}). In the following, we record some simple
facts about length spectra. Most of them are direct consequences of
the theorem. An in-depth treatment of this notion is given
in~\cite{lspec}.

Let $(k_i)_{i=1}^s$ be the list of odd factors of $n$ in ascending
order. Thus, $k_1=1$, $k_2$ (if it exists) is the smallest odd prime
factor of $n$, and $n = 2^dk_s$ for some $d \ge 0$. Let $r$ be the
largest index ($1 \le r \le s$) such that $k_r < \sqrt{2n}$.
\begin{enumerate}
\item The length of an even decomposition of $n$ is of the form
  $2n/k_j$, and hence the highest power of $2$ that divides any even
  element of $\lspec(n)$ is $2^{d+1}$. This observation rules out, for
  example, the possibility of the set $\{1,2,3,4\}$ being a length
  spectrum.
  
\item The smallest number with a length spectrum of size $s$ is
  $3^{s-1}$. The smallest number with $m$ in its length spectrum is
  clearly 
  \[ 1+2+\cdots + m = m(m+1)/2.\] In other words, $m \in \lspec(n)$
  implies $m(m+1)/2 \le n$. Hence $m \le (-1+\sqrt{1+8n})/2$. This
  gives an upper bound on the elements of $\lspec(n)$ in terms of $n$.
  
\item The longest decomposition of $n$ has length $\max\{k_r,
  2n/k_{r+1}\}$ (the maximum is $k_r$ if $n$ has no even
  decompositions). The shortest non-trivial decomposition of $n$ (if
  any) has length $\min\{k_2, 2n/k_s\} = \min\{k_2, 2^{d+1}\}$.
  
\item The condition $k_s < \sqrt{2n}$, or equivalently $k_s <
  2^{d+1}$, is clearly both necessary and sufficient for $n$ to have
  only odd decompositions. The situation, however, is quite different
  for even decompositions: if $k_j > \sqrt{2n}$, then $k_s/k_j <
  2k_s/k_j < \sqrt{2n}$, so the number of even decompositions of $n$
  is at most the number of odd decompositions of $n$. Consequently, if
  every non-trivial decomposition of $n$ is even, then $n$ has exactly
  one non-trivial decomposition. Those numbers with this property are
  precisely the numbers of the form $2^dk$ with $k$ an odd prime $>
  2^{d+1}$.
\end{enumerate}

\section*{Epilogue}
Finding the exact number of decompositions of $n$ can be hard. In
general, the formula $\prod_{p}(e_p+1)$ is impractical for large $n$
since it essentially calls for the prime factorization of $n$.  What
LeVeque obtained in~\cite{lev}, among other things, was the average
order of the number of decompositions as a function of $n$.  Moreover,
he discussed not only the sums of consecutive integers but the sums of
arithmetic progressions in general.  Readers with a taste for analytic
number theory will find his article enjoyable.

Guy gave a very short proof of the theorem in~\cite{guy}, then deduced
from it a characterization of primes. He gave some rough estimates of
the number of decompositions and also remarked that finding this
number explicitly is not easy.

Weil's book~\cite{weil} is merely a collection of exercises for the
elementary number theory course given by him at University of Chicago
in the summer of 1949. Maxwell Rosenlicht, an assistant of Weil's at
that time, was in charge of the ``laboratory'' section and responsible
for most of the exercises. It was a relief to learn from Weil that the
challenge of motivating students to work on problems is rather common.
The following is part of the foreword taken directly from that book:
\begin{quote}
  The course consisted of two lectures a week, supplemented by a
  weekly ``laboratory period'' where students were given exercises...
  . The idea was borrowed from the ``Praktikum'' of German
  universities.  Being alien to the local tradition, it did not work
  out as well as I had hoped, and student attendance at the problem
  sessions soon became desultory.
\end{quote}

\begin{figure}[htp]
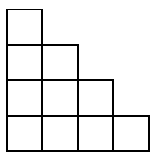

\begin{equation*}
  \yng(1,2,3,4)
\end{equation*}
\caption{\small The $4$-th triangular number is $10$.}
\label{f:tri}
\end{figure}

\begin{figure}[htp]
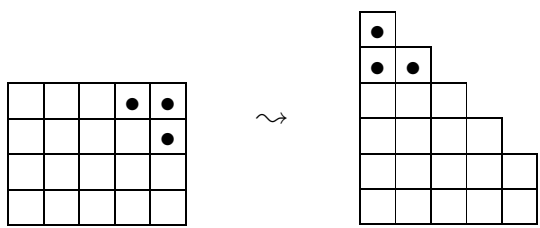

\begin{equation*}
  \raise-.52in
  \hbox{\young(\hfil\hfil\hfil\bullet\bullet,\hfil\hfil\hfil\hfil\bullet,\hfil\hfil\hfil\hfil\hfil,\hfil\hfil\hfil\hfil\hfil)}
  \qquad \leadsto \qquad
  \young(\bullet,\bullet\bullet,\hfil\hfil\hfil,\hfil\hfil\hfil\hfil,\hfil\hfil\hfil\hfil\hfil,\hfil\hfil\hfil\hfil\hfil)
\end{equation*}
\caption{\small Transforming a $4 \times 5$ rectangle into a trapezoid} 
\label{f:transf}
\end{figure}

\begin{figure}[htp]
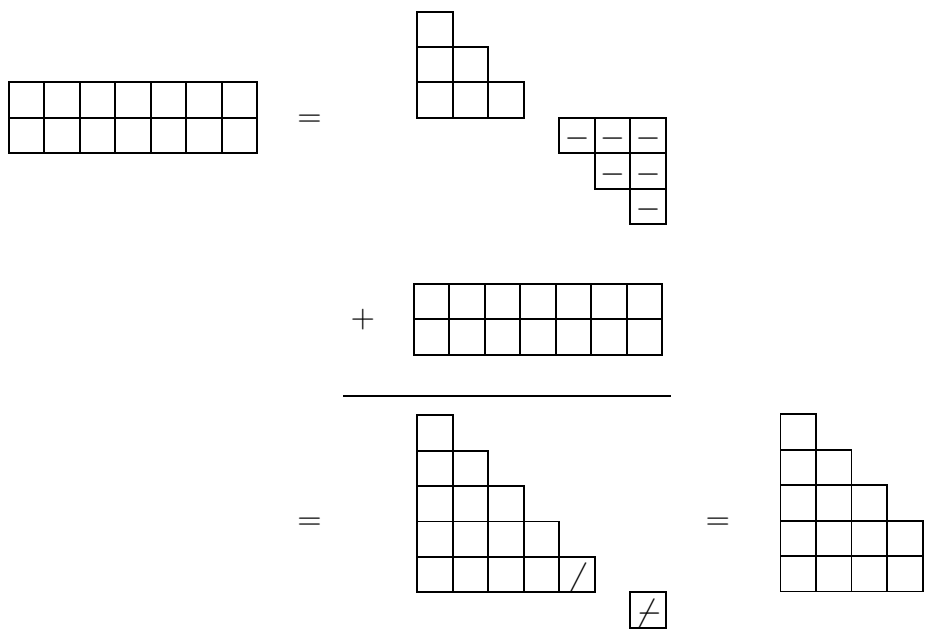

  \begin{align*}
    \yng(7,7) \quad &= \quad \phantom{+} \quad
    \young(\hfil,\hfil\hfil,\hfil\hfil\hfil,::::---,:::::--,::::::-)
    \quad \\ \\
    &\phantom{=} \quad + \quad \yng(7,7)\\
    &\phantom{=} \quad \underline{\phantom{111111111111111111111}} \\
    &= \quad \phantom{+} \quad 
    \young(\hfil,\hfil\hfil,\hfil\hfil\hfil,\hfil\hfil\hfil\hfil,\hfil\hfil\hfil\hfil\cancel,::::::\cneg)
    \quad = \quad \raise-23.5pt \hbox{
      \young(\hfil,\hfil\hfil,\hfil\hfil\hfil,\hfil\hfil\hfil\hfil,\hfil\hfil\hfil\hfil)}
  \end{align*}
  \caption{\small The decomposition (2,3,4,5) of 14.}
  \label{f:decomp}
\end{figure}

An obvious but crucial point in our proof of the conjecture is that
$0$ can be expressed as the sum of an odd number of consecutive
integers. Let me explain the intuition behind this trick here.  The
sum of the first $n$ consecutive natural numbers is called the {\em
  $n$-th triangular number} because $1, 2, \cdots, n$ can be
arranged to form a triangle (see Figure~\ref{f:tri}).
Now a sum of consecutive numbers can be viewed as the difference of
two triangular numbers (or a ``trapezoidal number'').  Also, a number
$n$, as long as it is not a power of 2, can be represented by an
$(n/k) \times k$ rectangle with $k>1$ an odd factor of $n$.  So the
question is: how can you get a trapezoid from such a rectangle?  Well,
it does not take a big leap of imagination to see that this can be
done by cutting off a corner of the rectangle and flipping it over.
For example, the diagrams in Figure~\ref{f:transf} illustrate how we
obtain the decomposition $(2,3,4,5,6)$ of $20$.
Of course, one has to worry about the case when $n/k < (k-1)/2$, but
this is exactly why the use of negative numbers comes in handy. The
diagrams in Figure~\ref{f:decomp} should be self-explanatory.
 
{\em Acknowledgment.}  I would like to thank Greg Kallo, the editor,
and the referees for many helpful suggestions on the presentation of
this article.

\end{document}